%% file: degreebound.tex
\documentclass{jssc}

\def\textsubscript#1%
{$_{\text{#1}}$}

\def\cdd{\mbox{\boldmath$\cdot$}~}
\newcommand{\rulex}{\hfill\rule{1mm}{3mm}}

\input{amssym.def}
\include{graphix}
\def\Z{{\mathbb{Z}}}
\def\bP{{\mathbb{P}}}
\def\tdeg{{\rm tdeg}}
\def\ord{{\rm ord}}
\def\bfa{{\bf a}}
\def\frakp{{\mathfrak p}}
\def\frakh{{\mathfrak h}}
\def\msindex{{\rm m.s.index}}
\def\overkt{{\overline{k(t)}}}

\input jsscN.tex
\begin{document}

%*************************************************************************************************************
% \biaoti{THE CAPITALIZED TITLE OF YOUR ARTICLE$^*$}{The list of authors' names with the LAST NAME capitalized
% and the authors' names should be separated by "\cdd"}{the first author's name \\ the first author's affiliation
% and Email address\\ the second author's name\\ the second author's affiliation. More can be listed like this.}
% {$^*$ The titles and numbers of the foundations that support this article.}
%*************************************************************************************************************
\title{Rational Solutions of First Order Algebraic Ordinary Differential Equations$^*$}%%%   Main Title of your paper  %%%
{\uppercase{Feng} Shuang \cdd \uppercase{Shen}
Li-Yong}%%% The names of the authors  %%%
{\uppercase{Feng} Shuang  \cdd \uppercase{Shen} Li-Yong\\
School of Mathematical Sciences, University of Chinese Academy of Sciences, Beijing $100049$, China\\
Email: fengshuang@ucas.ac.cn, lyshen@ucas.ac.cn    % Academy of Mathematics and Systems Science, Chinese Academy of Sciences, Beijing $100190$, China
} %%% The address of the authors  %%%
{$^*$This research was supported by Beijing Natural Science Foundation (Z190004), by National Key Research and Development
Project 2020YFA0713703, and by the Fundamental Research Funds for the Central Universities.\\
}

%*************************************************************************************************************
%The submission date of your article. For example: \drd{Received: June 8, 2006}
%*************************************************************************************************************

%*************************************************************************************************************
% The page header of the article.
% \dshm{Year}{Volume}{The capitalized RUNNING HEAD of your article with less than 48 letters}{The capitalized
% AUTHORS list with $\cdot$ separating different names or one can type "The name of the first author et al."
% if there are more than 4 authors.}
%*************************************************************************************************************

%*************************************************************************************************************
% \dab{The abstract}{Keywords}
%*************************************************************************************************************
%-------------------------------------------------------------------------
\Abstract{Let $f(t,y,y')=\sum_{i=0}^n a_i(t,y)y'^i=0$ be an irreducible first order ordinary differential equation with polynomial coefficients.
Eremenko in 1998 proved that there exists a constant $C$ such that every rational solution of $f(t,y,y')=0$ is of degree not greater than $C$.  Examples show that this degree bound $C$ depends not only on the degrees of $f$ in $t,y,y'$ but also on the coefficients of $f$ viewed as the polynomial in $t,y,y'$.  In this paper, we show that if $f$ satisfies $\deg(f,y)<\deg(f,y')$ or
$$
     \max_{i=0}^n \{ \deg(a_i,y)-2(n-i)\}>0
$$
then the degree bound $C$ only depends on the degrees of $f$ in $t,y,y'$, and furthermore we present an explicit expression for $C$ in terms of the degrees of $f$ in $t,y,y'$.}      % the abstract

        % the keywords

%\MRSubClass{05B05, 05B25, 20B25}      % MR(2000) Subject Classification

%\baselineskip 15pt

\section{Introduction}
The study of first order algebraic ordinary differential equations (AODEs in short) has a long history, which can be at least tracked back to the time of Fuchs and Poincar\'e.  Fuchs presented a sufficient and necessary condition so called Fuchs' criterion for a first order AODE having no movable singularity. Roughly speaking, an AODE is said to have movable singularities if it has a solution (with arbitrary constants) whose branch points depend on arbitrary constants. For instance the solution $y=\sqrt{t+c}$ of $2yy'-1=0$ has branch points $t=-c$, where $c$ is an arbitrary constant, so $2yy'-1=0$ has movable singularities.
Based on differential algebra developed by Ritt \cite{ritt} and the theory of algebraic function field of one variable, Matsuda \cite{matsuda} reproduced many classic results of first order AODEs. In particular, he presented an algebraic definition of movable singularities. In 1998, combining Matsuda's results and height estimates of points on a plane algebraic curve, Eremenko \cite{eremenko} showed that rational solutions of  first order AODEs have bounded degrees. In \cite{feng-feng}, we proved that if a first order AODE has movable singularities then it has only finitely many rational solutions. As for algebraic solutions of first order AODEs, Freitag and Moosa \cite{freitag-moosa} showed that  they are of bounded heights.

On the other hand, the algorithmic aspects of computing closed form solutions of AODEs have been extensively studied in the past decades. Several algorithms have been developed for computing closed form solutions (e.g. liouvillian solutions) of linear homogeneous differential equations (see \cite{barkatou,kovacic,vanhoeij-ragot-ulmer-weil,singer1,vanderput-singer} etc). Yet, the situation is different in the nonlinear case. Existing algorithms are only valid for AODEs of special types. Based on parametrization of algebraic curves, Aroca et al \cite{aroca-cano-feng-gao,feng-gao} gave two complete methods for finding rational and algebraic solutions of first order autonomous AODEs. Their methods were generalized by Winkler and his collegues to the class of first order non-autonomous AODEs whose rational general solutions involved arbitrary constants rationally as well as some other certain classes of AODEs (see \cite{vo-grasegger-winkler1,vo-grasegger-winkler2,chau-winkler, winkler} etc). Particularly, in \cite{vo-grasegger-winkler1}, the authors introduced a class of first order AODEs called maximally comparable AODEs and presented an algorithm to compute a degree bound for rational solutions of this kind of equations as well as first order quasi-linear AODEs.   Readers are referred to \cite{winkler} for a survey of recent developments in this direction.
Theoretically, it suffices to compute a degree bound for all rational solutions of  a first order AODE to find all its rational solutions by substituting an expression for $y$ with indeterminate coefficients and solving the resulting system of algebraic equations. The following example implies that the degrees of rational solutions may depend not only on the degrees of the original equation but also on its constant coefficients.
\begin{example}\label{ex:nomovable}
Let $m$ be an integer. Then $y=t^m$ is a rational solution of $ty'-my=0$.  The degree of $t^m$ depends on the constant coefficient $m$ of $ty'-my$.
\end{example}
%Notice that the differential equation in Example \ref{ex:nomovable} has no movable singularity, and then one sees that the degrees of rational solutions of such equations may rely on the constant coefficients. Our major concern in this paper is the first order AODEs having movable singularities, which occupy a large part of first order AODEs and only have finitely many rational solutions.

Let $f=\sum_{i=0}^n a_i(t,y)y'^i=0$ be an irreducible first order AODE with $n=\deg(f,y')$. Set
\begin{equation*}
\label{eqn:index}
     \msindex(f)=\max_{i=0}^n \{ \deg(a_i,y)-2(n-i)\}.
\end{equation*}
Fuchs' theorem (see Theorem 2 on page 11 of \cite{matsuda}) implies that $f=0$ has movable singularities if $\msindex(f)>0$. On the other hand, it was proved in \cite{eremenko} that if $f=0$ has movable singularities then it can be transformed into an AODE $g$ with positive $\msindex$. This motivates us to focus on first order AODEs with positive $\msindex$.
We prove that if $f$ satisfies $\deg(f,y)<\deg(f,y')$ or $\msindex(f)>0$, then the degrees of rational solutions of $f=0$ are independent of the constant coefficients of $f$. Furthermore using height estimates of points on a plane algebraic curve in \cite{feng-feng2}, we present an explicit degree bound in terms of the degrees of $f$ in $t,y,y'$.

The paper is organized as follows. In Section 2, we introduce some basic materials used in the later sections. In Section 3, we apply the results of height estimates in \cite{feng-feng2} to rational solutions of first order AODEs and obtain the main results in the paper. Furthermore, we consider some special types of first order AODEs where more compact bounds of rational solutions can be achieved.

Throughout this paper, $\Z$ stands for the ring of integers, $k$ for an algebraically closed field of characteristic zero, $k(t)$ for the field of rational functions in $t$ and $\overkt$ for the algebraic closure of $k(t)$. $\bP^m(\cdot)$ denotes the projective space of dimension $m$ over a field, and $(a_0:a_1:\dots:a_m)$ denotes a point in $\bP^m(\cdot)$ with coordinates $a_i$. As usual, for a polynomial $f(x_1,x_2,\dots,x_m)$, we use $\tdeg(f)$ and $\deg(f,x_i)$ to denote the total degree of $f$ and the degree of $f$ with respect to $x_i$ respectively.

\section{Basic materials}
In this section, we will introduce some basic materials used in this paper, including differential algebra, algebraic function fields of one variable and heights. Readers are referred to \cite{ritt, matsuda, chevalley, lang, serre} for details.

\subsection{Differential fields associated to AODEs}
All rings in this subsection are assumed to be commutative rings with unity.
\begin{definition}
A derivation on a ring $R$ is a map $\delta: R \rightarrow R$ satisfying that for all $a,b\in R$,
$$
   \delta(a+b)=\delta(a)+\delta(b),\,\,\delta(ab)=\delta(a)b+a\delta(b).
$$
A ring (resp. field) equipped with a derivation is called a differential ring (resp. differential field).  An ideal $I\subset R$ is called a differential ideal if $\delta(I)\subset I$. And a prime ideal of $R$ is called a prime differential ideal if it is a differential ideal.
\end{definition}
The field $k(t)$ can be endowed with a structure of differential field whose derivation $\delta$ is the usual derivation with respect to $t$, i.e. $\delta=\frac{\rm d}{{\rm d} t}$. Set $y_0=y$ and denote
$$k(t)\{y\}=k(t)[y_0,y_1,\dots]$$
 where $y_0,y_1,\dots$ are indeterminates.  One can extend the derivation $\delta$ on $k(t)$ to a derivation $\bar{\delta}$ on $k(t)\{y\}$ by assigning $y_i=\bar{\delta}^i(y_0)$ so that $k(t)\{y\}$ becomes a differential ring.  For the sake of notations, we use $\delta$ in place of $\bar{\delta}$. Elements in $k(t)\{y\}$ are called differential polynomials over $k(t)$. Let $f$ be a differential polynomial not in $k(t)$. Then there is a unique integer $d$ such that
$f\in k(t)[y_0,y_1,\dots,y_d]\setminus k(t)[y_0,y_1,\dots,y_{d-1}]$. This integer is called the order of $f$. We shall use $[\cdot]$ (resp. $\langle \cdot \rangle$) to stand for the differential (resp. algebraic) ideal generated by a set of differential polynomials (resp. polynomials) respectively.
When we say that $f$ is an irreducible differential polynomial, we mean that $f$ is irreducible over $k(t)$ as an algebraic polynomial. Suppose that $f$ is an irreducible differential polynomial.  Set
 $$
   \Sigma_f=\left\{ A\in k(t)\{y\} | \,\exists\, m>0\, \,\mbox{s.t.}\, S^m A^m \in [f] \right\}
$$
where $S=\partial f/\partial y_d$ and $d$ is the order of $f$. It was proved on page 30 of \cite{ritt} that $\Sigma_f$ is a prime differential ideal and so $k(t)\{y\}/\Sigma_f$ is a differential domain.
Lemma 2.2 of \cite{feng-feng} implies that the field of  fractions of $k(t)\{y\}/\Sigma_f$ is isomorphic to that of $k(t)[y_0,y_1,\dots,y_d]/\langle f \rangle$. Under this isomorphism, the field of fractions of $k(t)[y_0,y_1,\dots,y_d]/\langle f\rangle$ can be endowed with a structure of differential field. We shall still use $\delta$, or $'$ in short,  to denote the induced derivation on the field of fractions of $k(t)[y_0,y_1,\dots,y_d]/\langle f \rangle$.

In this paper, we mainly consider first order AODEs of the following form
\begin{equation*}
\label{eq:differentialeqn}
 f(y,y')=0
\end{equation*}
where $f(y,y')\in k(t)[y,y']\setminus k(t)$ is irreducible. As mentioned above, the field of fractions of $k(t)[y,y']/\langle f(y,y') \rangle$ is not only an algebraic function field over $k(t)$ but also a differential field, which is called a differential algebraic function field of one variable over $k(t)$ (see the definition on page 2 of \cite{matsuda}).
\begin{definition}
An element $r(t)\in k(t)$ satisfying
$f(r(t),r'(t))=0$ is called a rational solution of $f(y,y')=0$. The degree of $r(t)$, denoted by $\deg(r(t))$, is defined to be the maximum of the degrees of the denominator and numerator of $r(t)$.
\end{definition}

To avoid triviality, we always assume that the degrees of rational solutions which we consider in this paper are at least one. In next subsection we will introduce the concept of heights and one sees that the height of a rational function is exactly its degree.

\subsection{Heights in algebraic function fields of one variable}
First of all, we shall introduce some basic materials about algebraic function fields of one variable.

Let $L\subset \overkt$ be a finite extension of $k(t)$. Then $L$ is an algebraic function field of one variable over $k$. A  discrete valuation ring of $L$ over $k$ is a subring $V$ satisfying that
\begin{enumerate}
\item[1)] $k\subset V\neq L$; and
\item[2)] if $a\in L \setminus V$, then $a^{-1}\in V$.
\end{enumerate}
All non-invertible elements of $V$ form a maximal ideal $\frakp$ which is called a  place of $L$, and $V$  is called the corresponding ring of $\frakp$. Let $V$ be a discrete valuation ring with $\frakp$ as the place. There is an element $u\in V$, called a uniformizing variable of $\frakp$ or $V$, such that $\frakp=uV$ and $\bigcap_{n=1}^\infty u^n V=\{0\}$.

Let $\frakp$ be a place of $L$ and $V$ the corresponding ring of $\frakp$. Assume that $u$ is a uniformizing variable of $\frakp$. Then for every nonzero element $a$ of $L$, there is a unique integer $n$ such that $a=u^nv$ for some invertible element $v\in V$.
It is easy to see that the integer $n$ is independent of the choice of uniformizing variables.
Such $n$ is called the order of $a$ at $\frakp$ and denoted by $\ord_{\frakp}(a)$. We make the convention to write $\ord_{\frakp}(0)=\infty$. Then the place $\frakp$ induces a map $\ord_{\frakp}$ from $L \setminus \{0\}$ to $\Z$ sending $a$ to $\ord_{\frakp}(a)$. This map $\ord_{\frakp}$ is called the order function at $\frakp$.

\begin{remark}
In \cite{walker}, a place is presented by an equivalent class of irreducible parametrizations of the corresponding algebraic curve of $L$. Actually, these two definitions of places are equivalent to each other and the readers are referred to Remark 2.2 in \cite{feng-feng2} for details.
\end{remark}

Specially, $k(t)$ is the simplest algebraic function field of one variable over $k$. In the following example one can see what places, uniformizing variables and order functions will be when $L=k(t)$.

\begin{example}
\label{ex:algfunctionfield}
Let $L=k(t)$. When we express an element $a\in L$ as $g/h$ with $g,h\in k[t]$, we default that $\gcd(g,h)=1$. For each $c\in k$, set
$$
V_{c}=\left\{\frac{g}{h}\;\big |\; g,h\in k[t],\,h(c)\neq 0\right\}.
$$
It is clear that $k \subset V_{c} \subsetneq  L$. If $a\in L\setminus V_{c}$, there exist $g,h\in k[t]$ such that $a=g/h$ with $h(c)=0$, and thus $a^{-1}=h/g \in V_{c}$ since $g(c)\neq 0$.
Therefore $V_{c}$ is a discrete valuation ring of $L$. Similarly, one can verify that the set
$$
V_{1/t}=\left\{\frac{g}{h}\;\big |\; g,h\in k[t],\,\deg(h)\geq \deg(g)\right\}
$$
is also a discrete valuation ring of $L$. Actually, $V_{c}$ for any $c\in k$ and $V_{1/t}$ exhaust all the discrete valuation rings of $L$. And the corresponding places of them are respectively
$$
\frakp_{c}=\left\{\frac{g}{h}\;\big |\; g,h\in k[t],\,g(c)=0,\, h(c)\neq 0\right\}
$$
and
$$
\frakp_{1/t}=\left\{\frac{g}{h}\;\big |\; g,h\in k[t],\,\deg(h)>\deg(g)\right\},
$$
where a uniformizing variable of $\frakp_{c}$ (resp. $\frakp_{1/t}$) is $t-c$ (resp. $1/t$).

Let $a\in L\setminus \{0\}$. For any $c\in k$, $a$ can be written as
$$
a=(t-c)^m\frac{g_c}{h_c},
$$
where $m\in \Z,\,g_c,h_c\in k[t],\,g_c(c)h_c(c)\neq 0$, and therefore $\ord_{\frakp_{c}}(a)=m$. Assume that $a=g/h$ with $g,h\in k[t]$, then $\ord_{\frakp_{1/t}}(a)=\deg(h)-\deg(g)$.
\end{example}

Next, we introduce the definitions and properties of heights in algebraic function fields of one variable. Readers can refer to \cite{feng-feng2} for detailed information.

\begin{definition}
\label{def:height}
\begin{enumerate}
\item[1)]
Given $\bfa=(a_0:a_1: \dots :a_m)\in \bP^m(\overkt)$, let $L$ be a finite extension of $k(t)$ containing all $a_i$. We define the absolute logarithmic height (or simply height) of $\bfa$, denoted by $\frakh(\bfa)$,  to be
$$
 \frac{\sum_{\frakp} \max_{i=0}^m\{-\ord_{\frakp}(a_i)\}}{[L:k(t)]}
$$
where $\frakp$ ranges over all places of $L$.
\item[2)] For $a\in \overkt$, we define the height of $a$ to be $\frakh((1:a))$, denoted by $\frakh(a)$.
\item[3)] Let $f$ be a nonzero polynomial in $\overkt[x_1,x_2,\dots,x_m]$. We define the height of $f$ to be
\[
    \frakh(f)=\begin{cases} 0, & \mbox{$f$ contains exactly one term},\\
    \frakh(\bfa), & \mbox{otherwise}, \end{cases}
\]
where $\bfa$ is the point in some projective space whose coordinates are the coefficients of $f$.
\end{enumerate}
\end{definition}

Although there are infinitely many extensions $L$ of $k(t)$ containing the coordinates $a_i$, $\frakh(\bfa)$ is independent of the choices of $L$. Furthermore, $\frakh(\bfa)$ is independent of the choices of the homogeneous coordinates of $\bfa$. From Remark 2.6 in \cite{feng-feng2}, one sees that if $a\in \overkt$, then
$$
\frakh(a)=\frac{\deg(h,t)}{\deg(h,x)},
$$
where $h(t,x)$ is the nonzero irreducible polynomial over $k$ such that $h(t,a)=0$.
In particular, if $a\in k(t)$ then $\frakh(a)=\deg(a)$.

\begin{proposition}
\label{prop:heights2}
Suppose that $\bfa=(a_0:a_1: \dots :a_m)\in \bP^m(\overkt)$, $a_0,a_1,\dots,a_m\in k[t]$  and $\gcd(a_0,a_1,\dots,a_m)=1$. Then
$$
\frakh(\bfa)=\max_{i=0}^{m} \{\deg(a_i)\}.
$$
\end{proposition}
\proof
Let $L=k(t)$, then it is easy to see that
$$
\frakh(\bfa)=\sum_{\frakp} \max_{i=0}^m \{-\ord_{\frakp}(a_i)\}
$$
where $\frakp$ ranges over all the places of $L$. From Example \ref{ex:algfunctionfield}, one sees that every place of $L$ has a uniformizing variable of the form $1/t$ or $t-c$ for some $c\in k$.  Suppose that $\frakp$ is a place of $L$ with $t-c$ as a uniformizing variable. Then $\ord_\frakp(a_i)>0$ if and only if $t-c$ can divide $a_i$. Since $\gcd(a_0,a_1,\dots,a_m)=1$, there exists $i_0\in \{0,1,\dots,m\}$ such that $\ord_\frakp(a_{i_0})=0$, which implies that
$$\max_{i=0}^{m} \{-\ord_\frakp(a_i)\}=0.$$
Suppose that $\frakp$ is a place of $L$ with $1/t$ as a uniformizing variable, then
one has that $\ord_\frakp(a_i)=-\deg(a_i)$ and thus
$$\max_{i=0}^{m} \{-\ord_\frakp(a_i)\}=\max_{i=0}^{m}\{\deg(a_i)\}.$$
Consequently, $\frakh(\bfa)=\max_{i=0}^{m}\{\deg(a_i)\}$.
\rulex

\begin{example}
Compute the heights of following polynomials in $\overkt[x,y]$.
\begin{enumerate}
\item[1)] Let $f(x,y)=x+\sqrt{t}y$. Since the irreducible polynomial in $t,\sqrt{t}$ is $X^2-t$, one has that $$\frakh(f)=\frakh(\sqrt{t})=1/2.$$

\item[2)] Let $f(x,y)=(t^2+1)x/t^2+(t^3+1)y/t$. It follows from Proposition \ref{prop:heights2} that
$$\frakh(f)=\frakh(((t^2+1)/t^2:(t^3+1)/t))=\frakh((t^2+1:t^4+t))=4.$$
\end{enumerate}
\end{example}

The following two theorems are main results about height estimates of points on a plane algebraic curve in \cite{feng-feng2} (see Corollary 4.2 and Theorem 4.6 in \cite{feng-feng2}). The first one is a special case while the second is the general case.

\begin{theorem}\label{thm:boundpoints2}
Suppose $f(x,y)=\sum_{i=0}^m \sum_{j=0}^n c_{i,j}x^iy^j$ with $c_{i,j}\in \overkt,\,m=\deg(f,x),\, n=\deg(f,y)$. Assume that for all $0\leq i\leq m$ and $0\leq j\leq n$ if $c_{i,j}\neq 0$ then $mj+ni\leq mn$. Then for each $a,b\in \overkt$ with $f(a,b)=0$, one has that
$$
m\frakh(a)-mn\frakh(f)\leq n\frakh(b) \leq m\frakh(a)+mn\frakh(f).
$$
\end{theorem}

\begin{theorem}
\label{thm:boundpoints}
Let $f(x,y)$ be an irreducible polynomial in $\overkt[x,y]$ of degree $m$ with respect to $x$ and of degree $n$ with respect to $y$.  Suppose that $\rho=\tdeg(f)$ and $0<\epsilon<1$.  Then for every $a,b\in \overkt$ with $f(a,b)=0$, one has that
$$
  (1-\epsilon)m\frakh(a)-C\leq n\frakh(b)\leq (1+\epsilon)m\frakh(a)+C
$$
where
$$
C= 75\cdot 2^{13}\cdot (1/\epsilon)^6(\rho+1)^{\frac{40(\rho+1)^9}{\epsilon^3}} \frakh(f).
$$
\end{theorem}

In Section \ref{sec:mainresults}, we will use above two theorems to estimate degree bounds of rational solutions of first order AODEs.

\section{Main results}\label{sec:mainresults}
%Given a first order AODE $f(y,y')=0$ with coefficients in $k(t)$,
%viewing $y,y'$ as algebraic indeterminates the equation $f(y,y')=0$ can be considered as a plane algebraic curve over $\overkt$. If $f(r(t),r'(t))=0$ with $r(t)\in k(t)$, then $(r(t),r'(t))$ can be seen as a point on this curve. Using the results of height inequalities of points on a plane algebraic curve in \cite{feng-feng2}, we can estimate the degree bound of $r(t)$ by the total degree and the height of $f$.

Assume that $f(y,y')=0$ is an irreducible first order AODE with coefficients in $k(t)$ and that $r(t)\in k(t)$ is a rational solution of $f(y,y')=0$.
First, we shall show that if $\deg(f,y)> 2\deg(f,y')$ or $\deg(f,y)<\deg(f,y')$, then using Theorem \ref{thm:boundpoints} we can estimate a degree bound of $r(t)$ in terms of the total degree of $f$ in $y,y'$ and the height of $f$. Furthermore, with proper reduction $f$ can be assumed to be irreducible in three variables $t,y,y'$ and thus the height of $f$ is $\deg(f,t)$ by Proposition \ref{prop:heights2}. It means that $\deg(r(t))$ can be bounded only by the degrees of $f$ in $t,y,y'$.

The following lemma will play a key role in the proof of Theorems \ref{thm:boundforsols1} and \ref{thm:boundforsols2}.
Although they are existing results (see Corollary 2.10 in \cite{feng-feng2} and Lemma 3.9 in \cite{feng-gao}), we reprove these results for completeness.

\begin{lemma}\label{lem:rationalfunction}
Let $r(t)\in k(t)\setminus k$ and $c_1,c_2,c_3,c_4\in k$ with $c_1c_4-c_2c_3 \neq 0$. Then
\begin{enumerate}
  \item[$1)$] $\deg((c_1r+c_2)/(c_3r+c_4))=\deg(r)$;
  \item[$2)$] $\deg(r)-1\leq \deg(r')\leq 2\deg(r)$.
\end{enumerate}
\end{lemma}
\proof
Assume $r(t)=p(t)/q(t)$ with $p,q\in k[t],\,\gcd(p,q)=1$. Suppose $m=\deg(p(t))$ and $n=\deg(q(t))$.

1) Set
$$
\bar{r}=\frac{c_1r+c_2}{c_3r+c_4}=\frac{c_1p+c_2q}{c_3p+c_4q}.
$$
Since $\gcd(p,q)=1$ and $c_1c_4-c_2c_3 \neq 0$, it is clear that $\gcd(c_1p+c_2q,c_3p+c_4q)=1$. If $m\neq n$, it is easy to see $\deg(\bar{r})=\deg(r)$; if $m=n$, then at least one of the coefficients of $t^m$ in $c_1p+c_2q$ and $c_3p+c_4q$ is nonzero since $c_1c_4-c_2c_3 \neq 0$, and therefore $\deg(\bar{r})=\deg(r)$.

2) Since
$$r'(t)=\frac{p'(t)q(t)-p(t)q'(t)}{q^2(t)},$$
one sees that $\deg(r')\leq 2\deg(r)$. It is clear that $\deg(r')=\deg(r)-1$ if $n=0$. Next we assume that $n>0$. We shall prove that $\gcd(p'q-pq',q^2)=\gcd(q,q').$
Let $d=\gcd(q,q')$. Assume that $q=dq_1,\,q'=dq_2$, where $q_1,q_2\in k[t],\,\gcd(q_1,q_2)=1$ and $q_1$ is the square free part of $q$. Therefore one sees that $\gcd(p'q_1-pq_2,q^2)=1$
and thus
$$\gcd(p'q-pq',q^2)=\gcd(d(p'q_1-pq_2),q^2)=d.$$
If $m\leq n$, then $\deg(r')=\deg(q^2)-\deg(d)\geq 2n-n=n$, otherwise $\deg(r')=\deg(p'q-pq')-\deg(d)\geq m+n-1-n=m-1$ since $\deg(p'q-pq')=m+n-1$. In both cases, one sees that $\deg(r')\geq \deg(r)-1$.
\rulex

\begin{theorem}
\label{thm:boundforsols1}
Assume that $f(y,y')=0$ is an irreducible first order AODE with coefficients in $k(t)$, and assume further that $\deg(f,y)>2\deg(f,y')$ or $\deg(f,y)<\deg(f,y')$. If $r(t)$ is a rational solution of $f(y,y')=0$, then
$$
   \deg(r(t))\leq
   \begin{cases}
   75\cdot 2^{13}\cdot(\rho+1)^{40(\rho+1)^{12}+7}\frakh(f), \quad & \deg(f,y)>2\deg(f,y'),\\
   75\cdot 2^{13}\cdot(\rho+1)^{40(\rho+1)^{12}+7}\frakh(f)+(\rho+1)^2, \quad & \deg(f,y)<\deg(f,y'),
   \end{cases}
$$
where $\rho=\tdeg(f)$.
\end{theorem}
\proof
Denote $m=\deg(f,y),\,n=\deg(f,y')$.
Suppose that
$$
   f=h_1h_2\cdots h_s
$$
where $h_i$ is irreducible over $\overkt$. Since $f$ is irreducible over $k(t)$, one has that all $h_i$ are conjugate to each other and then
\begin{align*}
       \deg(h_i, y)=\deg(f,y)/s=m/s,\,\,\deg(h_i, y')=\deg(f,y')/s=n/s.
\end{align*}
Furthermore, $\frakh(h_i)\leq \frakh(f)$ due to Proposition 2.15 in  \cite{feng-feng2}. Assume that $r(t)$ is a rational solution of $f(y,y')=0$ then $r(t)$ is a rational solution of all $h_i=0$. Without loss of generality, assume $h_1(r(t),r'(t))=0$. Denote $\tilde{\rho}=\tdeg(h_1),\,\tilde{m}=\deg(h_1,y),\,\tilde{n}=\deg(h_1,y')$. Then
\begin{align*}
    \tilde{\rho} =\tdeg(f)/s=\rho/s,\,\tilde{m}=m/s,\,\tilde{n}=n/s.
\end{align*}

First, assume $m>2n$ and then $\tilde{m}>2\tilde{n}$.
Notice that $\frakh(r(t))=\deg(r(t))$. Set $\epsilon=1/(\tilde{\rho}+1)$.
Then using Theorem \ref{thm:boundpoints} and Lemma \ref{lem:rationalfunction} we have
$$
\frac{\tilde{\rho}}{\tilde{\rho}+1}\tilde{m} \deg(r(t))-C\leq \tilde{n}\deg(r'(t))\leq 2\tilde{n}\deg(r(t)),
$$
where
$
C\leq 75\cdot 2^{13}\cdot(\tilde{\rho}+1)^{40(\tilde{\rho}+1)^{12}+6}\frakh(f).
$
Therefore
\begin{equation*}
\label{eq:leftinequality}
    \left(\frac{\tilde{\rho}}{\tilde{\rho}+1}\tilde{m} -2\tilde{n} \right)\deg(r(t))\leq C,
\end{equation*}
where
\begin{align*}
    \frac{\tilde{\rho}}{\tilde{\rho}+1}\tilde{m} -2\tilde{n}=\frac{\tilde{\rho}(\tilde{m}-2\tilde{n})-2\tilde{n}}{\tilde{\rho}+1}\geq \frac{1}{\tilde{\rho}+1}
\end{align*}
since $\tilde{\rho}\geq \tilde{m} >2\tilde{n}$.
Then one sees that
\begin{align*}
    \deg(r(t))&\leq (\tilde{\rho}+1) C \\
    &\leq 75\cdot 2^{13}\cdot(\tilde{\rho}+1)^{40(\tilde{\rho}+1)^{12}+7}\frakh(f)\\
    &\leq 75\cdot 2^{13}\cdot(\rho+1)^{40(\rho+1)^{12}+7}\frakh(f).
    %&\leq 2^{20}\cdot(\rho+1)^{40(\rho+1)^{12}+7}\frakh(f).
\end{align*}

Next, assume $m<n$ and then $\tilde{m}<\tilde{n}$. Let $\epsilon,C$ be as above. Using Theorem \ref{thm:boundpoints} and Lemma \ref{lem:rationalfunction} again, we have
$$
\tilde{n}(\deg(r(t))-1)\leq \tilde{n}\deg(r'(t))\leq \frac{\tilde{\rho}+2}{\tilde{\rho}+1}\tilde{m} \deg(r(t))+C.
$$
A similar discussion one gets that
\begin{align*}
    \deg(r(t))&\leq (\tilde{\rho}+1)(C+\tilde{n})\\
    &\leq 75\cdot 2^{13}\cdot(\tilde{\rho}+1)^{40(\tilde{\rho}+1)^{12}+7}\frakh(f)+(\tilde{\rho}+1)^2\\
    &\leq 75\cdot 2^{13}\cdot(\rho+1)^{40(\rho+1)^{12}+7}\frakh(f)+(\rho+1)^2.
    %&\leq 2^{20}\cdot(\rho+1)^{40(\rho+1)^{12}+7}\frakh(f),
\end{align*} \rulex
%\begin{remark}
%Actually, from above proof process it can be seen that the condition $\frakh(f)>0$ in Theorem \ref{thm:boundforsols1} is only used in the case $\deg(f,y)<\deg(f,y')$. In other words, $\frakh(f)=0$ is allowed in another case $\deg(f,y)>2\deg(f,y')$.
%\end{remark}

Assume that
$
f(y,y')=\sum_{i=0}^n a_i(y)y'^i=0
$
is an irreducible first order AODE, where $a_i(y)\in k(t)[y]$ and $n=\deg(f,y')$.
Denote by
\begin{equation}\label{eq:msindex}
       \ell=\msindex(f)=\max_{i=0}^n \{ \deg(a_i)-2(n-i)\}.
\end{equation}
On the one hand, it is clear that $\ell>0$ if $\deg(f,y)>2\deg(f,y')$. On the other hand, we shall prove that if $\ell>0$ we can transform $f(y,y')=0$ into another irreducible first order AODE $g(z,z')=0$ with $\deg(g,z)>2\deg(g,z')$, and then use Theorem \ref{thm:boundforsols1} to estimate a degree bound of rational solutions of $f(y,y')=0$.
Suppose that $\ell>0$. Pick $c\in k$ such that $a_0(c)\neq 0$.
Set $y=(cz+1)/z$.  Then $y'=-z'/z^2$. Set
$$
     b_i(z)=a_i((cz+1)/z)z^{\ell+2n-2i} (-1)^i
$$
where $i=0,1,\dots,n$. Then an easy calculation yields that
\begin{align*}
   g(z,z')=\sum_{i=0}^n b_i(z)z'^i=z^{2n+\ell}f\left(\frac{cz+1}{z}, \frac{-z'}{z^2}\right).
\end{align*}
It is clear that $\deg(g,z')=n$, and $\deg(g,z)=\deg(b_0)=2n+\ell>2n$ since $a_0(c)\neq 0$.
Furthermore, $\tdeg(g)=2n+\ell$ because
$$
  2n+\ell\leq \tdeg(g)=\max_i\{ \deg(b_i)+i\} \leq \max_i\{2n+\ell-i\}=2n+\ell.
$$

Furthermore, we claim that $g(z,z')$ is irreducible over $k(t)$. First of all, assume that $\ell=\deg(a_{i_0})-2(n-i_0)$ for some $0
\leq i_0 \leq n$. Then we have that
$$
    b_{i_0}(0)=(-1)^{i_0}\cdot {\rm lc}(a_{i_0})\neq 0,
$$
where ${\rm lc}(a_{i_0})$ is the leading coefficient of $a_{i_0}$ with respect to $y$.
If $\gcd(b_0,b_1,\dots,b_n)\neq 1$, then the $b_i(z)$ have common zeroes and none of common zeroes is zero. It is easy to see that $(c\eta+1)/\eta$ is a common zero of all $a_i(y)$ if $\eta$ is a common zero of all $b_i(z)$. This contradicts with the fact that $\gcd(a_0,a_1,\dots,a_n)=1$.  Secondly, if $g(z,z')$ has a factor with positive degree in $z'$ then $f(y,y')$ will have a factor with positive degree in $y'$, a contradiction. This proves our claim.

Remark that $r(t)$ is a nontrivial rational solution of $g(z,z')=0$ if and only if $(cr(t)+1)/r(t)$ is a nontrivial rational solution of $f(y,y')=0$, and $\deg(r(t))=\deg((cr(t)+1)/r(t))$ by Lemma \ref{lem:rationalfunction}. That is to say, the degree bound of rational solutions of $f(y,y')=0$ is the same as that of $g(z,z')=0$.
Using Theorem \ref{thm:boundforsols1}, we will give an explicit degree bound of rational solutions of a first order AODE with positive $\msindex$.
\begin{theorem}
\label{thm:boundforsols2}
Assume that $f(y,y')=0$ is an irreducible first order AODE with coefficients in $k(t)$ and $\ell=\msindex(f)>0$. If $r(t)$ is a rational solution of $f(y,y')=0$, then
$$
   \deg(r(t))\leq (2\rho)^{(2\rho)^{15}}\frakh(f)
$$
where $\rho=\tdeg(f)$.
\end{theorem}
\proof
Notice that $\rho\geq 2$ since $f(y,y')=0$ has movable singularities, and if $\rho=1$ then $f(y,y')=0$ is a Riccati equation which has no movable singularity.
The notation being as above, from the previous discussion we only need to focus on the irreducible first order AODE $g(z,z')=0$. One sees that $\frakh(g)\leq \frakh(f),\, \deg(g,z')=n$ and
$
  \tdeg(g)=\deg(g,z)=2n+\ell>2n.
$
Then it follows from (\ref{eq:msindex}) that
$$
2n+\ell=\max_{i=0}^n \{ \deg(a_i)+2i\}\leq \max_{i=0}^n \{ \deg(a_i)+i\}+n \leq  2\rho-1,
$$
where the last inequality holds because $\ell>0$ induces $\rho>n$.
Assume that $r(t)$ is a rational solution of $f(y,y')=0$.
Due to Theorem \ref{thm:boundforsols1} we have
\begin{align*}
   \deg(r(t))&\leq 75\cdot 2^{13}\cdot(\tdeg(g)+1)^{40(\tdeg(g)+1)^{12}+7}\frakh(f)\\
             &\leq 2^{20}\cdot(2\rho)^{40(2\rho)^{12}+7}\frakh(f)\\
             &\leq (2\rho)^{(2\rho)^{15}}\frakh(f).
\end{align*}
\rulex

\begin{remark}
\begin{enumerate}
\item[1)] Multiplying the coefficients of $f$ by a common denominator and reducing them appropriately, we can assume that $f$ is an irreducible polynomial in three variables $t, y, y'$. Then it follows from Proposition \ref{prop:heights2} that $\frakh(f)=\deg(f,t)$. It means that, for differential polynomials $f$ in Theorems \ref{thm:boundforsols1} and \ref{thm:boundforsols2}, we can only use the degrees of $f$ in $t, y, y'$ to estimate a degree bound of rational solutions of $f=0$.

\item[2)] Theorem~\ref{thm:boundforsols2} implies that an autonomous first order AODE $f=0$ with positive $\msindex$ has no nontrival rational solutions, because $\frakh(f)=0$. In fact, suppose that $f=0$ has a nontrival rational solution. Then it will have infinitely many rational solutions. By Corollary 4.6 of \cite{feng-feng}, $f=0$ has no movable singularity. However, as $f=0$ has positive $\msindex$, Fuchs' theorem implies that $f=0$ has movable singularities, a contradiction.
\end{enumerate}
\end{remark}

In \cite{vo-grasegger-winkler1}, the authors developed two algorithms to compute rational solutions of maximally comparable first order AODEs and first order quasi-linear AODEs respectively. Let us first recall the definition of maximally comparable first order AODEs. Suppose that $f(y,y')=\sum_{i,j}c_{i,j}y^iy'^j$ is a first order differential polynomial over $k(t)$. Denote
$$
S(f)=\{(i,j)\in \Z^2 \mid c_{i,j}\neq 0\}.
$$
If there is $(i_0,j_0)\in S(f)$ satisfying that $i_0+j_0\geq i+j$ and $i_0+2j_0>i+2j$ for every $(i,j)\in S(f)\setminus (i_0,j_0)$, then $f$ is called maximally comparable. The following example shows that their algorithms can not completely deal with the cases we consider in Theorems \ref{thm:boundforsols1} and \ref{thm:boundforsols2}.

\begin{example}\label{ex:maxcomparable}
\begin{enumerate}
\item[1)] Let
$$
f(y,y')=\alpha_1(t)y'^3+\alpha_2(t)y^2 y'^2+\alpha_3(t)
$$
where $\alpha_i(t)\in k[t]\setminus \{0\}$.
Then $S(f)=\{(0,3),(2,2),(0,0)\}$. Since $2+2\geq 0+3$ but $2+2\cdot 2=0+2\cdot 3$, $f$ is not maximally comparable.
However, it is obvious that $\deg(f,y)<\deg(f,y')$.
\item[2)] Let
$$
f(y,y')=\alpha_1(t)y^\ell y'^n+\alpha_2(t)y^{2n+\ell}+\alpha_3(t)
$$
where $\ell,n \in \Z, \, \ell,n>0,\,\alpha_i(t)\in k[t]\setminus \{0\}$.
Then $S(f)=\{(\ell,n),(2n+\ell,0),(0,0)\}$. Since $2n+\ell+0\geq \ell+n$ but $2n+\ell+2\cdot 0=\ell+2\cdot n$, $f$ is not maximally comparable.
However, one sees that $\msindex(f)=\ell>0$. Furthermore, compared with Theorems \ref{thm:boundforsols1} and \ref{thm:boundforsols2} we can compute a more compact bound for this case in Example \ref{ex:degreebound}.
\end{enumerate}
\end{example}

Next we consider some special types of first order AODEs, in which the degree bounds of rational solutions can be extensively reduced compared with those in Theorems \ref{thm:boundforsols1} and \ref{thm:boundforsols2}.

\begin{proposition}\label{prop:boundforsols}
Let $f(y,y')=\sum_{i,j}c_{i,j}y^iy'^j \in k(t)[y,y']$ be a first order differential polynomial. Assume $\deg(f,y)=m,\,\deg(f,y')=n$ and $mj+ni\leq mn$ for $c_{i,j}\neq 0$. Furthermore, assume $m>2n$ or $m<n$.  If $r(t)\in k(t)$ is a rational solution of $f(y,y')=0$, then
$$
\deg(r(t))\leq
\begin{cases}
mn\frakh(f)/(m-2n),\quad &m>2n,\\
(mn\frakh(f)+n)/(n-m), \quad &m<n.
\end{cases}
$$
\end{proposition}
\proof
If $m>2n$, due to Theorem \ref{thm:boundpoints2} and Lemma \ref{lem:rationalfunction} one sees that
$$
m\frakh(r(t))-mn\frakh(f)\leq n\frakh(r'(t))\leq 2n\frakh(r(t))
$$
and then $\frakh(r(t))\leq mn\frakh(f)/(m-2n)$.
If $m<n$, using Theorem \ref{thm:boundpoints2} and Lemma \ref{lem:rationalfunction} again we have
$$
n(\frakh(r(t))-1)\leq n\frakh(r'(t))\leq m\frakh(r(t))+mn\frakh(f)
$$
and then $\frakh(r(t))\leq (mn\frakh(f)+n)/(n-m)$.
\rulex

\begin{proposition}\label{prop:boundforsols2}
Assume that $f(y,y')=\sum_{i=0}^n a_i(y)y'^i\in k(t)[y,y']$ is an irreducible first order differential polynomial with $\ell=\msindex(f)>0$. Assume further that $a_n(y)=\alpha y^\ell,\,\alpha \in k(t)\setminus \{0\},\,y^\ell \, |\, a_i(y),\,1\leq i\leq n-1$ and $a_0(0)\neq 0$. If $r(t)\in k(t)$ is a rational solution of $f(y,y')=0$, then
$$
\deg(r(t))\leq \frac{n(2n+\ell)}{\ell}\frakh(f).
$$
\end{proposition}
\proof
Set $y=1/z$, then $y'=-z'/z^2$. Let $$g(z,z')=f\left(\frac{1}{z},\frac{-z'}{z^2}\right)z^{2n+\ell}=\alpha(-z')^n+\sum_{i=0}^{n-1} b_i(z)z'^i,$$
where $b_i(z)=a_i(1/z)z^{2n+\ell-2i}(-1)^i$.
It is clear that $\deg(g,z')=n,\,\deg(g,z)=\deg(b_0(z))=2n+\ell$ since $a_0(0)\neq 0$. For $1\leq i\leq n-1$, assume that $\deg(b_i(z))= 2n+\ell-2i-m_i$, where $m_i\geq \ell$ since $y^\ell \, |\, a_i(y)$. Therefore one sees that for $0\leq i\leq n-1$
$$
\deg(b_i(z))n+i(2n+\ell)=(2n+\ell)n+\ell i-m_i n\leq (2n+\ell)n,
$$
where the last inequality holds since $\ell \leq m_i,\,i< n$. If $r(t)\in k(t)$ is a rational solution of $g(z,z')=0$, then Proposition \ref{prop:boundforsols} implies that
$$
\deg(r(t))\leq \frac{n(2n+\ell)}{\ell}\frakh(g).
$$
Notice that $r(t)$ is a nontrivial rational solution of $g(z,z')=0$ if and only if $1/r(t)$ is a nontrivial rational solution of $f(y,y')=0$ and $\frakh(g)=\frakh(f)$. Then our assertion is proved.
\rulex

\begin{example}
\label{ex:degreebound}
Let
$$
f(y,y')=\alpha_1(t)y^\ell y'^n+\alpha_2(t)y^{2n+\ell}+\alpha_3(t),
$$
where $\ell,n \in \Z, \, \ell,n>0,\,\alpha_i(t)\in k[t],\,\alpha_1(t)\alpha_3(t)\neq 0$ and $\gcd(\alpha_1,\alpha_2,\alpha_3)=1$. It is clear that $\msindex(f)=\ell>0$.
Then $\frakh(f)=\max_{i=0}^3 \{\deg(\alpha_i)\}$ due to Proposition \ref{prop:heights2}.
For each rational solution $r(t)$ of $f(y,y')=0$, it induces from Proposition \ref{prop:boundforsols2} that
$$\deg(r(t))\leq \frac{n(2n+\ell)\max_{i=0}^3 \{\deg(\alpha_i)\}}{\ell}.$$
The following are three specific examples of computing rational solutions of $f(y,y')=0$ when $\ell=n=1$.
\begin{enumerate}
\item[1)] Let $\alpha_1(t)=\alpha_2(t)=1,\,\alpha_3(t)=t$, then $\deg(r(t))\leq 3$ and using the method of undetermined coefficients for $r(t)$ one sees that there are no rational solutions of this case.
\item[2)] Let $\alpha_1(t)=t,\,\alpha_2(t)=\alpha_3(t)=1$, then $\deg(r(t))\leq 3$ and using the same way one sees that all rational solutions of this case are constants $c$ with $c^3+1=0$.
\item[3)] Let $\alpha_1(t)=1,\,\alpha_2(t)=0,\,\alpha_3(t)=-t$, then $\deg(r(t))\leq 3$ and using the same way one sees that all rational solutions of this case are $\pm t$.
\end{enumerate}
\end{example}

%%Please make sure that your given name is abbreviated as the first capital letter, such as Zhang X T, Tami T,...

\end{document}

%% file: jsscN.tex
\makeatletter
\def\@oddfoot{\hfill}
\newcount\shumeicount
\def\setshumei#1#2#3{%
  \shumeicount=\count0
  \def\@oddhead{%
    \raise-5pt\hbox to0pt{\vrule width\hsize height 0pt depth 0.4pt\hss}\relax
    \ifnum \shumeicount=\count0
      \raise-7pt\hbox to0pt{\vrule width\hsize height 0pt depth 0.4pt\hss}\relax
      #1
    \else
      \ifodd\count0
        #2
      \else
        #3
       \fi
     \fi
  }%
}
\makeatother
\makeatletter
\def\@oddfoot{\hfill}
\newcount\shujiaocount
\def\setshujiao{%
  \shujiaocount=\count0
  \def\@oddfoot{%
      \ifodd\count0
      \else
      \fi
  }%
}
\makeatother
\def\title#1#2#3#4{{
  \vspace*{0.3cm}
  \begin{flushleft} \Large\bf #1\end{flushleft}
  \vspace*{-0.2cm}
      \begin{flushleft}
      \bf #2
      \end{flushleft}
      \footnotetext{\hspace{-6mm} #3\\ #4}}}
\def\dshm#1#2#3#4
{\setshumei{J Syst Sci Complex (#1) #2:
{\thepage--\pageref{LastPage}}\hfill}
            {\hfill {\small #3}\hfill\hbox to0pt{\hss\thepage}}
            {\hbox to0pt{\thepage\hss}\hfill {\small #4}\hfill
            }
            \setshujiao}
\def\drd#1#2
{{\vskip 1cm\small \begin{flushleft}
 #1 \\
 #2\\
\copyright The Editorial Office of JSSC \&  Springer-Verlag GmbH
Germany 2018
\end{flushleft}}}
\def\tilde{\widetilde}
\def\bar{\overline}
\def\epsilon{\varepsilon}
\def\proof{\vspace{1mm}\indent {\it Proof}\quad}

%% file: degreebound.bbl
\begin{thebibliography}{199}
\bibitem{ritt}
Ritt J F, {\it Differential Algebra}, American Mathematical Society, Providence, Rhode Island, 1950.

\bibitem{matsuda}
Matsuda M, {\it First Order Algebraic Differential Equations}, Springer-Verlag, Berlin, 1980.

\bibitem{eremenko}
Eremenko A, Rational solutions of first-order differential equations, {\it Ann. Acad. Sci. Fenn. Math.}, 1998, {\bf 23}(1): 181--190.

\bibitem{feng-feng}
Feng S and Feng R, Descent of ordinary differential equations with rational general solutions, {\it Journal of Systems Science and Complexity}, 2020, {\bf 33}: 2114-2123.

\bibitem{freitag-moosa}
Freitag J and Moosa R, Finiteness theorems on hypersurfaces in partial differential-algebraic geometry, {\it Advances in Mathematics}, 2017, {\bf 314}: 726--755.

\bibitem{singer1}
Singer M F, Liouvillian solutions of {$n$}th order homogeneous linear differential equations, {\it American Journal of Mathematics}, 1981, {\bf 103}(4): 661--682.

\bibitem{kovacic}
Kovacic J J, An algorithm for solving second order linear homogeneous differential equations, {\it Journal of Symbolic Computation}, 1986, {\bf 2}(1): 3--43.

\bibitem{barkatou}
Barkatou M A, On rational solutions of systems of linear differential equations, {\it Journal of Symbolic Computation}, 1999, {\bf 28}(4-5): 547--567.

\bibitem{vanhoeij-ragot-ulmer-weil}
Van Hoeij M, Ragot J F, Ulmer F, and Weil J A, Liouvillian solutions of linear differential equations of order three and higher, {\it Journal of Symbolic Computation}, 1999, {\bf 28}(4-5): 589--609.

\bibitem{vanderput-singer}
Van der Put M and Singer M F, {\it Galois Theory of Linear Differential Equations}, Springer-Verlag, Berlin, 2003.


\bibitem{aroca-cano-feng-gao}
Aroca J M, Cano J, Feng R, and Gao X S, Algebraic general solutions of algebraic ordinary differential equations, {\it Proceedings of the $2005$ International Symposium on Symbolic and Algebraic Computation}, New York, 2005.

\bibitem{feng-gao}
Feng R and Gao X S, A polynomial time algorithm for finding rational general solutions of first order autonomous {ODE}s, {\it Journal of Symbolic Computation}, 2006, {\bf 41}(7): 739--762.

\bibitem{chau-winkler}
Ng\^{o} L X C, and Winkler F, Rational general solutions of first order non-autonomous parametrizable {ODE}s, {\it Journal of Symbolic Computation}, 2010, {\bf 45}(12): 1426--1441.

\bibitem{vo-grasegger-winkler1}
Vo T N, Grasegger G, and Winkler F, Computation of all rational solutions of first-order algebraic {ODE}s, {\it Advances in Applied Mathematics}, 2018, {\bf 98}: 1--24.

\bibitem{vo-grasegger-winkler2}
Vo T N, Grasegger G, and Winkler F, Deciding the existence of rational general solutions for first-order algebraic {ODE}s, {\it Journal of Symbolic Computation}, 2018, {\bf 87}: 127--139.

\bibitem{winkler}
Winkler F, The algebro-geometric method for solving algebraic differential equations---a survey, {\it Journal of Systems Science and Complexity}, 2019, {\bf 32}(1): 256--270.

\bibitem{feng-feng2}
Feng R, Feng S, and Shen L Y, Quasi-equivalence of heights in algebraic function fields of one variable, http://arxiv.org/abs/2111.13025.

\bibitem{chevalley}
Chevalley C, {\it Introduction to the Theory of Algebraic Functions of One Variable}, American Mathematical Society, Providence, Rhode Island, 1963.

\bibitem{lang}
Lang S, {\it Fundamentals of Diophantine Geometry}, Springer-Verlag, New York, 1983.

\bibitem{serre}
Serre J P, {\it Lectures on the Mordell-Weil Theorem}, Springer Fachmedien Wiesbaden GmbH, 1997.

\bibitem{walker}
Walker R J, {\it Algebraic Curves}, Princeton University Press, Princeton, 1950.





%\bibitem{hess}
%Hess F, Computing {R}iemann-{R}och spaces in algebraic function fields and related topics, {\it Journal of Symbolic Computation}, 2002, {\bf 33}(4): 425--445.

%\bibitem{fulton}
%Fulton W, {\it Algebraic Curves}, Addison-Wesley Publishing Company,
%Redwood City, CA, 1989.
%
%\bibitem{huang-ierardi}
%Huang M D and Ierardi D, Efficient algorithms for the {R}iemann-{R}och problem and for addition in the {J}acobian of a curve, {\it Journal of Symbolic Computation}, 1994, {\bf 18}(6): 519--539.

%\bibitem{sendra-winkler}
%Sendra J R and Winkler F, Tracing index of rational curve parametrizations, {\it Computer Aided Geometric Design}, 2001, {\bf 18}(8): 771--795.

\end{thebibliography}
